\documentclass{amsart}
\usepackage{amssymb}
\usepackage{graphicx}
\usepackage[cmtip,all]{xy}

\DeclareFontFamily{U}{wncy}{}
\DeclareFontShape{U}{wncy}{m}{n}{<->wncyr10}{}
\DeclareSymbolFont{mcy}{U}{wncy}{m}{n}
\DeclareMathSymbol{\Sh}{\mathord}{mcy}{"58} 

\newtheorem{theorem}{Theorem}[section]
\newtheorem{lemma}[theorem]{Lemma}

\newtheorem{corollary}[theorem]{Corollary}

\theoremstyle{definition}

\newtheorem{example}[theorem]{Example}

\newtheorem*{theorem*}{Theorem}
\theoremstyle{remark}
\newtheorem{remark}[theorem]{Remark}

\numberwithin{equation}{section}

\begin{document}

\title{Indivisibility of class numbers of imaginary quadratic fields}


\author{Olivia Beckwith}
\address{}
\curraddr{}
\email{}
\thanks{}

\date{}

\maketitle 
\begin{abstract}
We quantify a recent theorem of Wiles on class numbers of imaginary quadratic fields by proving an estimate for the number of negative fundamental discriminants down to $-X$ whose class numbers are indivisible by a given prime $\ell$ and whose imaginary quadratic fields satisfy any given set of local conditions. This estimate matches the best results in the direction of the Cohen-Lenstra heuristics for the number of imaginary quadratic fields with class number indivisible by a given prime. This general result is applied to study rank 0 twists of certain elliptic curves.
\end{abstract}

\section{Introduction}
Ideal class numbers of imaginary quadratic fields have been studied since Gauss, who conjectured that for any given $h$, there are only finitely many negative fundamental discriminants $D$ such that $h(D) = h$. The history of Gauss' Conjecture is rich. The conjecture was shown to be true by work of Heilbronn \cite{Heilbronn}, who did not show how to find the imaginary quadratic fields with a given class number.  Siegel \cite{Siegel} proved that $h(-D)$ grows like $|D|^{1/2}$, but did so ineffectively. In other words, for each $\epsilon > 0$ he proved that for sufficiently large $D$ that there are positive  constants $c_1$ and $c_2$ for which
$$
c_1 D^{1/2 -\epsilon} < h(-D) < c_2 D^{1/2+\epsilon}
$$
While explicit upper bounds for $h(-D)$ are known, the constants $c_1$ are ineffective for all $\epsilon$. Baker \cite{Baker} and Heegner \cite{Heegner} computed the complete finite list of negative fundamental discriminants D for which $h(D) = 1$.  The works Gross and Zagier \cite{GrZa} and Goldfeld \cite{Goldfeld2} produce a lower bound for $h(D)$ which is asymptotically smaller than Siegel's bound, but is effective and allows one (in principle) to compute the complete list of imaginary quadratic fields with any given class number. 

It is natural to ask what else can be said about the structure of ideal class groups. For example, how often should we expect the $\ell$-torsion subgroup of the class group to be trivial for a given odd prime $\ell$? The {\it Cohen-Lenstra heuristics} \cite{CL} predict an answer: 
\begin{equation}\label{eq:CL}
\lim_{X \to \infty} \frac{ \# \{ - X < D < 0: \ell \nmid h(D) \}}{X} = \prod_{n = 1}^{\infty} \left( 1 - \frac{1}{\ell^{n}} \right) = 1 - \frac{1}{\ell} - \frac{1}{\ell^2} + \frac{1}{\ell^5} \cdots
\end{equation}
Here the numbers $D$ are fundamental discriminants. Note that the Cohen-Lenstra heuristics actually predict much more about the structure of the class groups, give similar predictions for real quadratic fields, and have been generalized by others to other number fields. For a concise description for the quadratic number field case, the reader is encouraged to read Chapter 5 Section 10 of \cite{Cohen}. 

Numerical data provides some evidence for the Cohen-Lenstra heuristics, and for $\ell = 3$, strong theorems supporting equation (\ref{eq:CL}) are known. Gauss' {\it genus theory} says that the number of order 2 elements of the class group is $2^{t-1} - 1$, where $t$ is the number of distinct prime divisors of the discriminant (see Proposition 3.11 of \cite{Cox}). For $\ell = 3$, a theorem of Davenport and Heibronn \cite{DH} says that if $\epsilon >0$, then for $X$ sufficiently large we have
$$
\frac{ \# \{-X < D < 0: 3 \nmid h(D) \}}{X} \ge \frac{1}{2} - \epsilon.
$$
They proved this by showing that the cubic number fields are in a discriminant preserving correspondence with a certain set of classes of binary cubic forms, and they used this fact to count the order $3$ elements of class groups of quadratic number fields.

For $\ell > 3$ much less is known about the $\ell$-torsion of class groups. Soundararajan \cite{Sound} used analytic techniques to count $\ell$-torsion points of class groups, and showed 
$$
\# \{ -X < D < 0 : \ell | h(D) \} \gg X^{\frac{1}{2} + \epsilon(\ell)},
$$
where $\epsilon (\ell) >0$ approaches 0 as $\ell \to \infty$.  Kohnen and Ono \cite{KO} used the theory of modular forms to study the occurrence of class groups with trivial $\ell$-torsion for $\ell > 3$. They proved for any $\epsilon >0$, for sufficiently large $X$ we have
$$
\# \{-X < D < 0 : \ell \nmid h(D) \} \ge \left( \frac{2 (\ell -2)}{\sqrt{3} (\ell - 1)} - \epsilon \right) \frac{\sqrt{X}}{\log{X}}.
$$

Information about the structure of class groups of quadratic fields can be used to study questions about Mordell-Weil groups of elliptic curves in families of quadratic twists, however, additional information about the splitting and ramification data of the quadratic number fields is often required for such applications. For $E: y^2 = p(x)$ an elliptic curve over $\mathbb{Q}$ with $p(x)$ in Weierstrass form, we define the twist of $E$ by a fundamental discriminant $D$ to be the elliptic curve defined by 
$$
E_D : y^2 D = p(x).
$$
Note that $E_D$ is isomorphic to $E$ over $\mathbb{Q} (\sqrt{D})$, but not over $\mathbb{Q}$. The {\it Heegner hypotheses} are a set of conditions about how the rational primes of bad reduction of an elliptic curve split in an imaginary quadratic field. The work of Kolyvagin  on the Birch and Swinnerton-Dyer Conjecture (see \cite{Ko1}, \cite{Ko2}) is based on the existence of suitable quadratic twists of elliptic curves in which the twisting discriminant satisfy prescribed Heegner hypotheses. Combining his work with an important theorem of Gross and Zagier, who showed that the height of the Heegner point is a multiple of the derivative of the $L$-series of the elliptic curve at 1, it follows that the Birch and Swinnerton-Dyer Conjecture holds when the analytic rank is at most 1.

Heegner points have played an important role in studying Goldfeld's Conjecture, which concerns the ranks of the twists as $D$ varies over the set of fundamental discriminants. Define $M^r_E (X) := \# \{D: |D|<X: \mbox{ord}_{s=1} L(s, E_D) = r \}$. If $E/\mathbb{Q}$ is an elliptic curve and $r$ is 0 or 1, then
$$
M^r_E (X) \sim \frac{X}{2}, \quad X \to \infty.
$$ 
The best general results on Goldfeld's Conjecture were, until recently, due to Perelli, Pomykala, and Skinner (see \cite{OS} and \cite{PePo}). For the rank 0 case, Ono and Skinner \cite{OS} showed that 
\begin{equation}\label{eq:OS} 
M^0_E (X) \gg \frac{ X}{\log {X}}.
\end{equation} 
For the rank 1 case, Perelli and Pomykala \cite{PePo} showed 
\begin{equation}\label{eq:PePo} M^1_E (X) \gg_{\epsilon} X^{1 - \epsilon} \end{equation}
for any $\epsilon > 0$. 

Earlier this year, Kriz and Li \cite{KL} showed for a large class of elliptic curves,
$$
M^r (X) \gg \frac{X}{\log^{\frac{5}{6}}(X)}.
$$

Strong results on Goldfeld's conjecture have been obtained for special elliptic curves by making use of the aformentioned theorem of Davenport and Heilbronn on the 3-indivisibility of class numbers. Using the half-integral weight modular forms established by Waldspurger and a theorem of Frey \cite{Frey}, James \cite{James} showed that the elliptic curve with Cremona label 14B satisfies $M^0_E (X) \gg X$. 

Showing that an elliptic curve has a positive proportion of twists with rank one requires more than Waldspurger's modular forms. Vatsal \cite{Vatsal} used a theorem of Gross and Zagier \cite{GrZa} to show that the elliptic curve $E = X_0(19)$ has $M^r_E (X) \gg X$ for $r = 0,1$.  Vatsal's argument was extended by Byeon \cite{Byeon} to elliptic curves in the isogeny class of an elliptic curve with a nontrivial cuspidal 3-torsion point and square-free conductor. 

The results toward Goldfeld's conjecture described above apply to certain elliptic curves with residually reducible mod 3 Galois representations, and rely on a refinement of the theorem of Davenport and Heilbronn due to Horie and Nakagawa \cite{HN}. Their refinement showed that a positive proportion of imaginary quadratic fields have trivial $\ell$-torsion and satisfy prescribed local conditions.  
One might hope to extend the work of Horie and Nakagawa to a theorem on $\ell$-indivisibility of class groups for $\ell > 3$ by refining the work of Kohnen and Ono \cite{KO} in an analogous way. 

A barrier to refining Kohnen and Ono's theorem is showing that the modular forms arising in their argument have Fourier coefficients which are supported on prescribed arithmetic progressions and are nontrivial modulo $\ell$. This is difficult because for many modular forms this property doesn't hold. For example, the values of the partition function $p(n)$ are the Fourier coefficients for the modular form $1 / \eta(z)$, and the Ramanujan congruences tell us $p(5n + 4) \equiv 0 \pmod{5}$, and so sieving the Fourier expansion of this form can return a modular form which is trivial modulo 5. Here $\eta(z):=q^{1/24}\prod_{n=1}^{\infty}(1-q^n)$ (here $q:=e^{2\pi i z}$ throughout) is Dedekind's eta-function, and it is a weight -1/2 weakly holomorphic modular form.

Recently, Wiles \cite{Wiles} established the existence of imaginary quadratic fields with prescribed local data whose class numbers are indivisible by a given odd prime $\ell$.

\begin{theorem*}{(Wiles)}
Let $\ell \ge 5$ be prime, and let $ S_0 , S_+ , S_{-}$ be finite disjoint sets of distinct odd primes not containing $\ell$ such that the following are true:

\begin{enumerate}
\item $S_0$ does not contain any primes which are $ 1 \pmod{\ell}$

\item $S_+$ does not contain any primes which are $ -1 \pmod{\ell}$

\item $S_-$ does not contain any primes which are $1 \pmod{\ell}$ and $-1 \pmod{4}$. 
\end{enumerate}
Then there exists a negative fundamental discriminant $D$ such that $\ell \nmid h(D)$, and $\mathbb{Q} (\sqrt{D})$ splits at every prime in $S_+$, is inert at every prime in $S_-$, and ramifies at every prime in $S_0$.
\end{theorem*}

In view of the work of Horie and Nakagawa when $\ell=3$ \cite{HN}, the goal of the present work is to prove a quantified version of the theorem of Wiles for the $\ell > 3$ case by obtaining an estimate for the number of imaginary quadratic fields which satisfy the conclusion of Wiles' theorem, similar to the estimate of Kohnen and Ono.

We define
\begin{equation}\label{eq:sturmbound}
M_{\Sigma} := \frac{1}{8} [ \Gamma_0 (1) : \Gamma_0 ( N_{\Sigma})]
\end{equation}
 and 
\begin{equation}\label{eq:level}
N_{\Sigma} : = 4 Q_{\Sigma}^6 (\prod_{q \in S_0 \cup S_- \cup S_+} q^6),
\end{equation}
where $Q_{\Sigma}$ is equal to $1$ if $S_-$ is nonempty and otherwise is the smallest odd prime not contained in $S_+ \cup S_- \cup S_0$ which is not congruent to 1 modulo $\ell$ and -1 modulo 4.

Our main theorem is the following estimate for the smallest discriminant divisible by a given prime $p$ lying in a certain arithmetic progression which satisfies the conclusion of Wiles' theorem.

\begin{theorem}\label{thm:maintheorem}
Suppose $p >  M_{\Sigma}$ is a prime such that the following are true:
\begin{enumerate}
\item We have that $\left( \frac{p}{\ell} \right) = 1$ and $p \not\equiv 1 \pmod{\ell}$,

\item We have that $p \equiv 1 \pmod{8}$,

\item For odd primes $q \le M_{\Sigma}$, $q \neq \ell$, we have $\left( \frac{p}{q} \right) = 1.$

\end{enumerate}
Then there is some $k_p \le p M_{\Sigma} $ such that $p \nmid k_p$ and $\ell \nmid h( - k_p p)$ and $\mathbb{Q} ( \sqrt{ - k_p p})$ ramifies at all  primes of $S_0$, splits at every prime in $S_+$, and is inert at every prime in $S_-$. 
\end{theorem}

Combining this result with Dirichlet's theorem on primes in arithmetic progressions, we obtain the following corollary, which can viewed as an extension of \cite{KO} to allow for local conditions. To state it, we let $T_{\Sigma, \ell}$ denote the set of all fundamental discriminants which satisfy the conclusions of Theorem \ref{thm:maintheorem}. That is, $T_{\Sigma, \ell}$ contains the set of negative fundamental discriminants $D$ of quadratic fields $K$ which ramify at all  primes of $S_0$, split at every prime in $S_+$, and are inert at every prime in $S_-$, and have $\ell \nmid h(D)$. Also, let $r_{\Sigma}$ be the number of odd primes less than $M_{\Sigma}$, excluding $\ell$. Then we have the following:
\begin{corollary}\label{thm:maincor}
Let $\ell$ be an odd prime.  If $\epsilon >0$, then for sufficiently large $X$ we have
$$
\# \{ -X < D <0: \ell \nmid h(D), D \in T_{\Sigma, \ell} \} \ge \left(  \frac{\ell -2 }{ (\ell - 1) 2^{r_{\Sigma}+4}  \sqrt{M_{\Sigma}}} - \epsilon \right) \frac{\sqrt{X}}{\log{X}}.
$$
\end{corollary}

One can apply Corollary \ref{thm:maincor} to count twists of elliptic curves which have Mordell-Weil rank 0 over $\mathbb{Q}$ and trivial $\ell$-Selmer group. To state this result, it is convenient to define the following subsets of primes dividing the conductor $N_E$. 
Let $\tilde{S}_E$ be the subset of odd primes dividing the conductor $N_E$ of $E$ defined by 
\begin{equation}\label{eq:SE}
\widetilde{S}_E := \{ p|N_E: p \equiv -1 \pmod{\ell}, \ell \nmid ord_{p}(\Delta_E) \},
\end{equation}
where $\Delta_E$ is the discriminant of $E$.
Also, we set
\begin{equation}\label{eq:Tplus}
T_+ = \{ p|N_E, ord_p(j_E) < 0; E/\mathbb{Q}_p \text{ is not a Tate Curve} \},
\end{equation}
and
\begin{equation}\label{eq:Tminus}
T_- = \{ p|N_E: p \notin T_+, p \equiv 3 \pmod{4} \}.
\end{equation}
A Tate curve $E / \mathbb{Q}_p$ is such that $E / \overline{\mathbb{Q}}_p \simeq \overline{\mathbb{Q}}_p^* / q^{\mathbb{Z}}$ for some $q \in \mathbb{Q}_p$, for details see Appendix C of \cite{Silverman}.

\begin{corollary}\label{thm:rank0}
Suppose $E/\mathbb{Q}$ is an elliptic curve with odd conductor $N_E$, and suppose  $E$ has a $\mathbb{Q}$-rational torsion point $P$ of odd prime order $\ell$, and suppose $P$ is not contained in the kernel of reduction modulo $\ell$.  Assume $ord_{\ell} (j(E)) \ge 0$. Also assume $\widetilde{S}_E = \emptyset$ and neither $T_+$ nor $T_-$ contain a prime which is $1 \pmod{\ell}$.

Then we have
$$
\# \{ -X < D <0: rk(E_D) = 0, Sel_{\ell} (E_D) = \{1\} \} \gg  \frac{\sqrt{X}}{\log{X}}.
$$
\end{corollary}

\begin{remark}
As mentioned above, the best general results on Goldfeld's Conjecture are due to Ono, Perelli, Pomykala, and Skinner (see equations \ref{eq:OS} and \ref{eq:PePo}). The corollary given here falls short of improving on this estimate. However, it is a refinement in that it gives rank 0 twists whose $\ell$-Selmer groups have trivial $\ell$-parts. This is the best known estimate for this type of problem. Corollary \ref{thm:rank0} presumably can be extended to also give rank 1 twists which simultaneously have trivial $\ell$-Shafarevich-Tate groups. This claim would require a careful study of the aforementioned paper of Frey \cite{Frey}.
\end{remark}

This paper is organized as follows. In Section 2, we describe a theorem of Zagier relating class numbers of imaginary quadratic fields to the coefficients of a weight $\frac{3}{2}$ mock modular form, and we use his result to prove a lemma which is vital to the proof of our main result. In Section 3, we prove Theorem \ref{thm:maintheorem} and Corollaries \ref{thm:maincor} - \ref{thm:rank0}, and in Section 4, we give examples to illustrate our results. 

\section*{Acknowledgements}
The author thanks Ken Ono for suggesting this project, and thanks Edray Goins and the referee for their many helpful comments.

\section{Hurwitz Mock Modular forms}
\subsection{Zagier's Eisenstein Series}
Throughout, $\mathbb{H}$ is the upper half plane, $z = x + iy$ is a complex number in $\mathbb{H}$ with $x,y \in \mathbb{R}$, and $q := e^{2 \pi i z}$. Also, for $ N \in \mathbb{Z}_{> 0}$, $k \in \frac{1}{2} \mathbb{Z}$, and $\chi$ a Dirichlet character, we let $M_{k} (\Gamma_0 (N), \chi)$ and $S_{k} (\Gamma_0 (N), \chi)$ denote the usual vector spaces of integer and half-integer weight modular forms and cusp forms.  

The proof of Theorem 1.1 requires generalizations of modular forms, the so-called {\it harmonic Maass forms}. We will describe only briefly the main properties of harmonic Maass forms. To learn more, the reader is encouraged to read \cite{BFOR} and \cite{Ono}. A harmonic Maass form is a real-analytic function which transforms like a modular form. All harmonic Maass forms have a natural decomposition, as
$$
f ( z) = f^+ ( z) + \frac{(4 \pi y)^{1-k}}{k-1} \overline{c_f (0) } + f^- ( z),
$$
where $f^+$ and $f^-$ have Fourier expansions as follows, for some $m_0 \in \mathbb{Z}$:
$$
f^+ ( z) = \sum_{n = m_0 }^{\infty} c_f^+ (n) q^n,
$$
and
$$
f^- (z) = \sum_{\substack{n > 0 }} \overline{c_f^- (n)}  \Gamma ( 1-k, 4 \pi n y) q^{-n},
$$
where $\Gamma_(s,x) := \int_x^{\infty} t^{s-1} e^{-t} dt$ is the \textit{incomplete Gamma-function}. 
The form $f^+$ is called the \emph{holomorphic part} of $f$, and $\frac{(4 \pi y)^{1-k}}{k-1} \overline{c_f(0) } + f^- ( z)$ is called the \emph{nonholomorphic part} of $f$. If the nonholomorphic part of $f$ is trivial, then $f$ is a weakly holomorphic modular form. When the nonholomorphic part is nontrivial, $f^+$ is called a \emph{mock modular form}. We let $H_k(\Gamma_0(N), \chi)$ denote the space of harmonic Maass forms with Nebentypus character $\chi$ on $\Gamma_0(N)$.

While many harmonic Maass forms have poles at cusps, not all of them do.  A harmonic Maass form $f \in H_k (\Gamma_0(N), \chi)$ is said to be of {\it moderate growth} if there exists $\epsilon >0$ such that
$$
f(z) = O(e^{\epsilon y})
$$
as $y \to \infty$, and if analogous conditions hold at all cusps of $\Gamma_0(N)$. We let $H_k^{mg} (\Gamma_0(N), \chi)$ denote the space of weight $k$ harmonic Maass forms of moderate growth. Harmonic Maass forms of moderate growth do not have poles at cusps. Moreover, harmonic Maass forms of moderate growth which have trivial nonholomorphic part are holomorphic modular forms. 

Throughout we let $D$ be a negative fundamental discriminant, and we let $h(D)$ be the class number for the quadratic field $\mathbb{Q}  ( \sqrt{D})$.  We use the Hurwitz class numbers $H(n)$, which are defined as follows. Suppose $-n = Df^2$, where $D < 0$ is a fundamental discriminant.
\begin{equation}\label{eq:hurwitz}
H(n) = \frac{h (D)}{w(D)} \sum_{d | f} \mu(d) \left( \frac{D}{d} \right) \sigma_1 \left(\frac{f}{d} \right),
\end{equation}
where $\sigma_1$ is the usual sum of divisors function and $w(-n)$ is half the number of units in the integer ring of $\mathbb{Q} (\sqrt{-n})$. 
%
 %
Let $\mathcal{H} (z)$ be defined by

\begin{equation}\label{eq:zagier}
 \mathcal{H}(z) := - \frac{1}{12} + \sum_{n=1}^{\infty} H(n) q^n + \frac{1}{8 \sqrt{\pi}} \sum_{n \in \mathbb{Z}}  \Gamma \left( - \frac{1}{2}, 4 \pi n^2 y \right) q^{-n^2},
\end{equation}
where $\Gamma(\alpha,x)$ is the usual incomplete Gamma-function. Zagier showed that $\mathcal{H}(z)$ is a harmonic Maass form \cite{Zagier}. In particular, if $\xi_{\frac{3}{2}}$ is the differential operator defined by $\xi_{\frac{3}{2}} := 2 i y^{\frac{3}{2}} \frac{\overline{\partial}}{\partial \overline{z}}$, then Zagier showed the following: 

\begin{theorem}{(Zagier)}\label{thm:zagier}
$\mathcal{H}(z)$ is a weight $\frac{3}{2}$ harmonic Maass form of moderate growth on $\Gamma_0(4)$. Moreover, $\xi_{3/2} (\mathcal{H}) = - \frac{1}{16 \pi} \Theta$, where $\Theta(z) := \sum_{n \in \mathbb{Z}} q^{n^2}$ is the Jacobi theta function. 
\end{theorem}

We use $\mathcal{H}(z)$ to  construct modular forms whose coefficients represent the fundamental discriminants which correspond to fields with the desired splitting conditions. Then we argue as in \cite{JO} and \cite{KO}.

\begin{remark}
The weight $3/2$ modular form $\sum_{n =0}^{\infty} r(n) q^n := \Theta(z)^3$ is intimately tied to class numbers for imaginary quadratic fields. It is well known that the $r(n)$ are given by Hurwitz class numbers $H(n)$. 
\begin{theorem}[Gauss]
$$
r(n) = \begin{cases}
12 H(n) \\
24 H(n) \\
r(n/4) & n \equiv 0 \pmod{4} \\
0 & n \equiv 7 \pmod{8}
\end{cases}
$$
\end{theorem}
The modular form $\theta^3$ was used in many previous results on indivisibility of class numbers (see for example \cite{KO}). However, it is insufficient for our result, because its Fourier coefficients are not supported on all arithmetic progressions. For the square free $n$ with $ n \equiv 7 \pmod{8}$, the class numbers $h (-n)$ are not represented. 
\end{remark}

\subsection{Sieving Zagier's Mock Modular Form}
We require the following result, which shows that we can define holomorphic modular forms whose coefficients are supported on fundamental discriminants satisfying local conditions and are given by class numbers. Given sets $S_+, S_-, S_0$ as in Theorem \ref{thm:maintheorem}, we let $A_{\Sigma}$ be defined as the set of positive integers $n$ such that the following hold:
\begin{enumerate}
\item For $p \in S_+ \cup S_- \cup S_-$, $p^2 \nmid n$.

\item $\mathbb{Q} (\sqrt{-n})$ splits at the primes in $S_+$, ramifies at the primes in $S_0$, and is inert at the primes in $S_-$.

\end{enumerate}

\begin{lemma}\label{thm:mainlemma}
Let $S_+, S_-, S_0$ be sets as in Theorem \ref{thm:maintheorem}, and assume that $S_-$ is nonempty.

Then there is a weight $\frac{3}{2}$ modular form $H^{\Sigma} (z) = \sum_{n=1}^{\infty} a(n) q^n$ on $\Gamma_0( N_{\Sigma})$, where $N_{\sigma}$ is as in equation \ref{eq:level}, such that
$$
a(n) = \begin{cases}
H(n) & n \in A_{\Sigma} \\
0 & \text{otherwise}
\end{cases}
$$
\end{lemma}

The idea is to take combinations of twists of Zagier's function $\mathcal{H}(z)$ to obtain holomorphic modular form. The key properties of $\mathcal{H}(z)$ that allow us to do this are
\begin{enumerate}
\item the Fourier expansion of the non-holomorphic part is supported on terms of the form $q^{-n^2}$, which allows us to use twisting to annihilate the non-holomorphic part of $\mathcal{H}(z)$, and 
\item $\mathcal{H} (z)$ has moderate growth at poles, which ensures that any linear combination of twists of $\mathcal{H}(z)$ will not have any exponential singularities, as a weakly holomorphic modular form would.
\end{enumerate}

For $\chi$ a Dirichlet character modulo $m$, the twist of $\mathcal{G}(z) := \sum_{n=0}^{\infty} a(n,y) q^n \in H_k (\Gamma_0(N), \psi)$ by $\chi$ is given by
$$
\mathcal{G}_{\chi} (z) = \sum_{n \in \mathbb{Z}} \chi(n) a(n,y) q^n .
$$
If $d$ is a positive integer, the operators $U(d), V(d)$ are defined, as one does when working with holomorphic modular forms, by 
$$
(\mathcal{G}| U(d)) (z) := \sum_{n \in \mathbb{Z}} a(dn,y) q^n
$$
and
$$
(\mathcal{G}| V(d)) (z) := \sum_{n \in \mathbb{Z}}^{\infty} a(n,y) q^{dn}.
$$

It is well known that a twist of a modular form is itself a modular form, for a proof, see Proposition 17 on page 127 of \cite{Koblitz}. The same proof shows that twists of harmonic Maass forms are also harmonic Maass forms. Specifically, for $\mathcal{G}(z) \in H^{mg}_{k+ \frac{1}{2}} ( \Gamma_0(4N), \psi)$, the form $\mathcal{G}_{\chi} (z)$ belongs to $H_{k + \frac{1}{2}}^{mg} ( \Gamma_0 (4N m^2), \psi \cdot \chi^2)$,  and $(\mathcal{G} | V(d) )$ and $(\mathcal{G}| U(d))$ lie in $H_{k+\frac{1}{2}}^{mg} ( \Gamma_0(4Nd), \psi \cdot \left( \frac{4d}{\cdot} \right))$.

\noindent \textit{Proof of Lemma \ref{thm:mainlemma}:}
First, we take a combination of twists for which the nonholomorphic part of $\mathcal{H} (z)$ is annihilated. Let $p$ be in $S_-$.  We have 
\begin{align*}
f(z) &:= \frac{1}{2} (\mathcal{H}(z) -  \left( \frac{-1}{p} \right) \mathcal{H}_{\left( \frac{\cdot}{p} \right)} (z)) \\
&= \sum_{n=1}^{\infty} \frac{1}{2} \left( 1 - \left( \frac{-n}{p} \right) \right) H(n) q^n + \frac{1}{16 \sqrt{\pi}} \sum_{p | n}  \Gamma(- \frac{1}{2}, 4 \pi n^2 y) q^{ - n^2}.
\end{align*}

Note that the coefficient of $q^n$ in the holomorphic part of $f(z)$is $\frac{1}{2} H(n)$ if $p | n$, $H(n)$ if $ (\frac{-n}{p} ) = -1$, and $0$ if $ \left( \frac{-n}{p} \right) = 1$. The nonholomorphic part of $f$ is supported on multiples of $p$, because twisting the nonholomorphic part by the Legendre symbol annihilates those coefficients. To eliminate what remains of the nonholomorphic part and the multiples of $p$ in the holomorphic part, we take the twist $ f_{\left( \frac{\cdot}{p}\right)^2}$. 

Repeating the above steps for every $p \in S_+ \cup S_-$, we obtain a form which is supported on $n$ for which the primes in $S_+ \cup S_-$ split or are inert in $\mathbb{Q} (\sqrt{-n})$ as desired. 

To obtain a modular form which is supported on coefficients which are multiples of the primes in $S_0$, let $d$ be the product of the primes in $S_0$. We apply the $U(d)$, operator, then twist by $\left( \frac{ -n}{q} \right)^2$ for each $q \in S_0$, and then apply the $V(d)$ operator.

\qed

\section{Proofs of Theorem \ref{thm:maintheorem} and Corollaries}
\subsection{Proof of Theorem \ref{thm:maintheorem}}
The proof of Theorem \ref{thm:maintheorem} requires a well known result of Sturm \cite{Sturm}, which says that if a modular form with integer Fourier coefficients is nonvanishing modulo a prime $\ell$, then there is a bound on the index of the first coefficient which is nonzero modulo $\ell$. To state his theorem, for a rational prime $\ell$ and a modular form $f(z) = \sum_{n=0}^{\infty} a(n) q^n \in M_{k} (\Gamma_0(N), \chi)$ with coefficients in $\mathbb{Z}$, we define $$ord_{\ell} (f):= \textit{min}_{n} \{ n: \ell \nmid a(n) \},$$
 and we say $ord_{\ell}(f) := \infty$ if $\ell | a(n)$ for all $n$.
\begin{theorem}[Sturm]\label{thm:sturm}
For a modular form $f(z) = \sum_{n=0}^{\infty} a(n) q^n \in M_{k} (\Gamma_0(N), \chi)$ with integer Fourier coefficients, if
$$
\textit{ord}_{\ell}(f) > \frac{k}{12} [ \Gamma_0(1):\Gamma_0(N)],
$$
then $\text{ord}_{\ell}(f) = \infty$.
\end{theorem}
\begin{remark}
Note that Sturm's theorem was originally only formulated for holomorphic modular forms of integer weight, but the proof carries over to half-integral weight modular forms.
\end{remark}
\noindent \textit{Proof of Theorem \ref{thm:maintheorem}:}
Let $H^{\Sigma} (z)$ be the modular form from Lemma \ref{thm:mainlemma} for $S_+$, $S_-$, and $S_0$, replacing $S-$ with $\{ Q_{\Sigma} \}$ if $S_- = \emptyset$. Let 
$$\mathcal{F} (z)  := \left(H^{\Sigma}| U(p) \right) -p \left( H^{\Sigma} | V(p) \right).$$

By Lemma \ref{thm:mainlemma}, the form $\mathcal{F}(z)$ is a modular form of weight $\frac{3}{2}$ on $\Gamma_0(p N_{\Sigma})$.  By Theorem 1 of Wiles \cite{Wiles}, $\mathcal{F} (z)$ has a Fourier coefficient which is indivisible by $\ell$. Therefore Sturm's Theorem tells us that we have
$$n_p := ord_{\ell} (\mathcal{F}) \le \frac{3}{24}  [ \Gamma_0(1) : \Gamma_0 (pN_{\Sigma} )] .$$ 

It follows from a well-known formula for $[\Gamma_0(1):\Gamma_0 (N)]$ (see for example \cite{CBMS}) that we have
$$n_p \le M_{\Sigma} (p+1).$$ 

We have that $n_p $ must be of the form $ f_p^2  k_p$, with $ k_p $ square free. It follows from conditions (1)-(3) in Theorem \ref{thm:maintheorem} that for all $n \le  M_{\Sigma}$, the $np^{th}$ Fourier coefficient of $\mathcal{F}(z)$ is divisible by $\ell$, so $p \nmid k_p$. Therefore either $- p k_p$ or $- 4 p k_p$ is a fundamental discriminant for an imaginary quadratic field satisfying the desired local conditions and whose class number is indivisible by $\ell$. 
\qed
\subsection{Proof of Corollary \ref{thm:maincor}}
Note that at least half of the values $k_p p$ from the main theorem must be distinct as $p$ varies over the primes greater than $M_{\Sigma}$ satisfying the conditions of Theorem \ref{thm:maintheorem}. If instead we had $k_p p = k_q q = k_r r$ with $p < q < r$, we would have $qr | k_p$, which would violate the bound on $k_p$. 

To count the fundamental discriminants down to $-X$ which satisfy the desired conditions, it suffices to count the primes which satisfy the conditions of Theorem \ref{thm:maintheorem} for which the fundamental discriminant from Theorem \ref{thm:maintheorem} is greater than $-X$.  

The primes $p$ that satisfy the third condition of Theorem \ref{thm:maintheorem} are those for which for each $q$ up to $M_{\Sigma}$, $p$ lies one of $\frac{q-1}{2}$ arithmetic progressions modulo $q$, which correspond to $p$ being a quadratic residue modulo $q$. Similarly, the other two conditions amount to restricting $p$ to certain arithmetic progressions modulo $2$ and $\ell$. 

For the fundamental discriminant corresponding to $p$ obtained from Theorem \ref{thm:maintheorem} to be guaranteed to be greater than $-X$, it suffices to require
$$
4 p M_{\Sigma} (p+1) \le X.
$$

It follows from Dirichlet's theorem for primes in arithmetic progressions that given $\epsilon > 0$  for sufficiently large $X$, we have
$$
\# \{ -X < D <0: \ell \nmid h(D), D \in T_{\Sigma} \} \ge \left( \frac{\ell -2}{ \ell -1} \frac{1}{ 2^{r_{\Sigma} + 4 }  \sqrt{M_{\Sigma}}} - \epsilon \right) \frac{\sqrt{X}}{\log{X}}.
$$

\subsection{Proof of Corollary \ref{thm:rank0}}
First we recall a theorem of Frey \cite{Frey}.
\begin{theorem}[Frey]\label{thm:frey}
Suppose $E / \mathbb{Q}$ is an elliptic curve with a $\mathbb{Q}$-rational torsion point $P$ of odd prime order $\ell$, and suppose $P$ is not contained in the kernel of reduction modulo $\ell$.  Suppose $\widetilde{S}_E = \emptyset$. Suppose that $D $ is a negative square-free integer coprime to $\ell N_E$ and satisfies
\begin{enumerate}
\item If $2 | N_E$ then $d \equiv 3 \pmod{4}$

\item If $ord_{\ell}(j(E))<0$, then $\left( \frac{D}{\ell} \right) = -1$,

\item If $p|N_E$ is an odd prime, then
$$
\left( \frac{d}{p} \right) = \begin{cases}
-1 & \textit{if } ord_p(j_E) \ge 0 \\
-1 & \textit{if } ord_p (j_E) < 0 \text{ and } E/\mathbb{Q}_p \text{ is a Tate curve} \\
1 & otherwise 
\end{cases}
$$
\end{enumerate}
Then $Sel_{\ell}(E_D)$ is nontrivial if and only if $\ell | h(D)$.
\end{theorem}

Now to prove the corollary, note that the twists $E_D$ have trivial $\ell$ torsion over $\mathbb{Q}$. We set
$$
S_+ = \{ p|N_E, ord_p(j_E) < 0: E/\mathbb{Q_p} \text{ is not a Tate Curve} \},
$$
and
$$
S_- = \{ p|N_E: p \notin S_+ \},
$$
and $S_0 = \emptyset$. It follows from Corollary \ref{thm:maincor} that there are at least $O(\frac{\sqrt{X}}{\log{X}})$ fundamental discriminants down to $-X$ which satisfy Frey's conditions, and so the result follows from Theorem \ref{thm:frey}.

\section{Examples}
Here we illustrate Theorem \ref{thm:maintheorem} and Corollary \ref{thm:rank0}. 
\begin{example}
Suppose that $\ell = 5$ and that the sets are $S_+ = \{ 3 \}$, $S_- = S_0 = \emptyset$. 

The smallest prime which satisfies the conditions of Theorem \ref{thm:maintheorem} is 394969. The smallest discriminant bounded by Theorem \ref{thm:maintheorem} is a multiple of this prime, however, it is clear that one shouldn't need to look at numbers that large to find imaginary quadratic fields which split at 3 and have a class number which is not divisible  by 5. By direct calculation, we see that for the primes $p$ less than 100, for all but 79 we have $5 \nmid h(-p)$, out of which 11 of the 21 corresponding imaginary quadratic fields split at 3. This discrepancy between the bounds predicted by Theorem \ref{thm:maintheorem} and the actual fundamental discriminants we observe is typical of these theorems, and it illustrates the main obstacles which remain in attacking the original Cohen-Lenstra conjectures.
\end{example}

\begin{example}
Let $E:y^2 + y = x^3 - x^2 + 20x -8$ be the elliptic curve with Cremona label 203.a1. Then $E(\mathbb{Q}) \simeq \mathbb{Z}/5\mathbb{Z}$. The conductor of $E$ is $7 \cdot 29$. It follows from Corollary \ref{thm:rank0} that we have
$$
\# \{ -X < D <0: rk(E_D) = 0, Sel_{5} (E_D) = \{1\} \} \gg  \frac{\sqrt{X}}{\log{X}}.
$$
\end{example}

%
%


\end{document}